# DS*: Tighter Lifting-Free Convex Relaxations for Quadratic Matching Problems

Florian Bernard[1,2]    Christian Theobalt[1,2]    Michael Moeller[3]

[1]MPI Informatics    [2]Saarland Informatics Campus    [3]University of Siegen

## Abstract

*In this work we study convex relaxations of quadratic optimisation problems over permutation matrices. While existing semidefinite programming approaches can achieve remarkably tight relaxations, they have the strong disadvantage that they lift the original $n{\times}n$-dimensional variable to an $n^2{\times}n^2$-dimensional variable, which limits their practical applicability. In contrast, here we present a lifting-free convex relaxation that is provably at least as tight as existing (lifting-free) convex relaxations. We demonstrate experimentally that our approach is superior to existing convex and non-convex methods for various problems, including image arrangement and multi-graph matching.*

## 1. Introduction

Matching problems that seek for correspondences between images, shapes, meshes or graphs are a long-standing challenge in computer vision and computer graphics. Computationally, they can be phrased as optimisation problems over binary variables that encode the matching. Whilst the formulation as a *discrete* optimisation problem appears most natural, in many scenarios a continuous formulation may be advantageous (*e.g.* for improving computational efficiency, or for representing uncertainties [46]).

In this work, we focus on convex relaxations of quadratic programming problems over permutation matrices, where we are particularly interested in finding both *scalable* as well as *tight* convex relaxations. To be more precise, we consider problems of the general form

$$\min_{X \in \mathbb{P}_n \cap \mathcal{C}} f(x) := x^T W x + c^T x, \quad (1)$$

where $x := \text{vec}(X) \in \mathbb{R}^{n^2}$ is the vector containing the columns of matrix $X \in \mathbb{R}^{n \times n}$, $\mathcal{C}$ is a (closed) convex set, and $\mathbb{P}_n$ is the set of $n{\times}n$ permutation matrices. Since

**Acknowledgements:** We thank A. Tewari for providing MoFA facial expression parameters [49]. This work was funded by the ERC Starting Grant CapReal (335545).

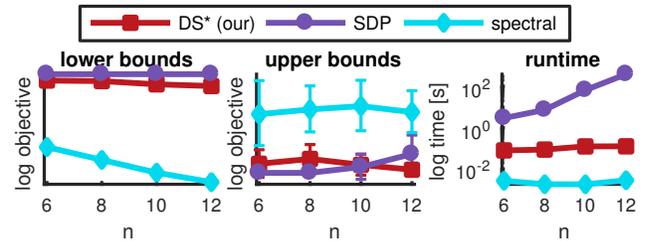

Figure 1. Normalised bounds (log-scale) and runtime in seconds (log-scale) for random instances of Problem (1). The (lifted) SDP has the best lower bounds but is not scalable. The spectral relaxation is efficient but has the weakest bounds. Our DS* method has reasonably good bounds and scales much better than the SDP.

in general problems of this form are known to be NP-hard [38], for moderately-sized problems one cannot expect to find a globally optimal solution. Thus, a lot of effort has been put into finding good solutions that may be suboptimal. Among them are semidefinite programming (SDP) relaxations of matching problems [59, 44, 27], which are known to produce good solutions for various binary problems [45, 53] by lifting the $n^2$-dimensional variable $x$ to an $n^4$-dimensional variable. Whilst SDP relaxations allow finding a solution in polynomial time (e.g. with roughly $\mathcal{O}(n^6)$ per-iteration complexity in SDCut [52]), the quadratic increase of the number of variables prohibits scalability (Fig. 1). The lifting-free methods include spectral [31] and convex relaxations [20, 1, 19, 17]. While they are better scalable, they achieve weaker bounds and thus usually do not result in high quality solutions. Our aim is to improve upon the tightness of lifting-free convex relaxations.

### 1.1. Related work

In this section we summarise existing works that are most relevant to our approach.

**Assignment problems:** The linear assignment problem (LAP) seeks to find a permutation matrix $X \in \mathbb{P}_n$ that minimises the (linear) objective $f_{\text{LAP}}(X) = \text{tr}(C^T X)$, where $C_{ij}$ indicates the cost of assigning object $i$ to object $j$ [35, 10]. The LAP is among the combinatorial methods that permit finding a global optimum in polynomial time,

*e.g.* using the Hungarian/Kuhn-Munkres [35] or the Auction algorithm [7]. However, a shortcoming of the LAP is that it neglects higher-order relationships. In contrast, the quadratic assignment problem (QAP) additionally takes the cost of matching pairs of objects (*e.g.* edges in a graph) into account. The Koopmans-Beckmann (KB) form [28] of the QAP refers to the the minimisation of

$$f_{\text{QAP}}(X) = \text{tr}(AXBX^T) + \text{tr}(C^T X) \qquad (2)$$

over permutation matrices $X \in \mathbb{P}_n$. Here, in addition to the linear costs encoded by $C$, the matrices $A$ and $B$ encode pairwise costs. In contrast to many existing works, we focus on the strictly more general Problem (1), which contains Problem (2) as special case (by setting $\mathcal{C} = \mathbb{R}^{n \times n}, c = \text{vec}(C)$, and $W = B^T \otimes A$, with $\otimes$ denoting the Kronecker product). As illustrated in [60], the KB form is limited to *scalar edge features* and *linear edge similarities*.

**Relaxation methods:** The QAP has received a lot of attention over decades (see *e.g.* [29, 10, 34]), which may be (at least partially) owed to the fact that it is NP-hard and that even finding an approximate solution within some constant factor of the optimal solution is only possible if P=NP [43]. Among the existing approaches that tackle the QAP are branch and bound methods [4] which rely on the inexpensive computation of bounds. In order to obtain such bounds, many relaxation methods have been proposed, most of which are either lifting-based and thus not scalable [59, 44, 27], or leverage the special structure of the Koopmans-Beckmann form [23, 38, 16, 2, 39, 40, 15, 36, 18] and are thus not directly applicable to the general form in Problem (1). A summary on various relaxations can be found in the survey paper by Loiola et al. [34]. Works that are both lifting-free and consider the general objective as in Problem (1) include spectral approaches [31, 14] and convex relaxation approaches [20, 17]. Our approach fits into the latter category, with our main contribution being the achievement of a relaxation that is provably at least as tight as the so-far tightest lifting-free convex relaxation [17].

**Graph matching:** The problem of bringing nodes and edges of graphs into correspondence is known as graph matching (GM). There are several GM variations, such as multi-graph matching [54], higher-order graph matching [30], or second-order graph matching, where the latter is an instance of the QAP [60]. Whilst there exist many different ways to tackle GM problems (*e.g.* [50, 12, 58, 26, 30, 47, 48]), in the following we focus on convex-to-concave path-following (PF) approaches, as they are most relevant in our context. The idea of PF methods is to approximate the NP-hard GM problem by solving a sequence of continuous optimisation problems. The motivation is to realise the rounding for obtaining binary solutions *within* the optimisation procedure, rather than as post-processing step. To this end, PF methods start with a (convex) subproblem where a (usually non-binary) globally optimal solution can be found. Then, using the current solution as initialisation for the next subproblem, they gradually move on to more difficult (non-convex) subproblems, until a binary solution is obtained when solving the last (concave) subproblem. The *PATH* algorithm [57] implements PF using convex and concave relaxations for the KB form of the QAP, as in (2). A similar approach is pursued in the *factorised graph matching* (FGM) method [60], where however the convex and concave relaxations are based on a factorised representation of the pairwise matching matrix $W$ in (1). In contrast to these methods, our approach directly establishes a convex relaxation without requiring any particular structure or factorisation of the pairwise matching matrix $W$.

### 1.2. Contributions

Our main contributions can be summarised follows:
1. We present a *novel convex relaxation framework* for quadratic problems over permutation matrices that *generalises* existing (lifting-free) relaxations and is *provably at least as tight* (Prop. 7).
2. To this end we provide a class of parametrised energy functions that are equivalent on the set of permutations (Prop. 5), infinitely many of them being convex.
3. Most importantly, we propose a proximal subgradient descent type algorithm to *efficiently find parameters that yield improved convex relaxations* (Alg. 1).
4. Our experimental validation confirms the flexibility and benefits of our method on various matching problems.

## 2. Background

In this section we present our notation followed by some background of convex relaxations.

### 2.1. Notation

We use $x := \text{vec}(X)$ to denote the vector containing all columns of the matrix $X$, and $X = \text{vec}^{-1}(x)$ reshapes the vector $x$ back into a matrix. $\mathbf{I}_n$ denotes the $n \times n$ identity matrix. $\mathbf{1}_n$ and $\mathbf{0}_n$ denote the constant $n$-dimensional vectors comprising of ones and zeros, respectively (subscripts may be omitted when the dimensionality is clear from context). We write $[n] := \{1, \ldots, n\}$ for $n \in \mathbb{N}$. Let $\mathbb{P}_n$ denote the set of $n \times n$ permutation matrices, *i.e.*

$$\mathbb{P}_n = \{X \in \{0,1\}^{n \times n} : X^T X = \mathbf{I}_n\}, \quad \text{and let}$$
$$\mathbb{DS}_n = \{X : X \geq 0, X\mathbf{1}_n = \mathbf{1}_n, \mathbf{1}_n^T X = \mathbf{1}_n^T\} \subset \mathbb{R}^{n \times n}$$

denote the set of doubly-stochastic matrices. For a symmetric matrix $W$ we use $W_-$ to denote the Euclidean projection onto the cone of negative semi-definite (NSD) matrices. The null space and the range of a matrix $A$ are denoted by $\ker(A)$ and $\text{im}(A)$, respectively. For a set $C$, by $\delta_C$ we

denote the indicator function that is given by

$$\delta_C(x) = \begin{cases} 0 & \text{if } x \in C, \\ \infty & \text{otherwise}. \end{cases} \quad (3)$$

## 2.2. Convex Relaxations

Convex relaxation methods are an important class of techniques for minimising non-convex energies. They aim at replacing the original cost function by a convex underapproximation whose minimiser ideally coincides with the global minimiser of the non-convex function, see Fig. 2 for an illustration. While the largest convex underapproximation, the *biconjugate* $F^{**}$, inherits such a property under weak additional conditions (e.g. by combining [24, Thm. 1.3.5] and the Krein-Milman theorem), it usually is at least as difficult to compute as solving the original problem to global optimality. One therefore has to settle for a smaller convex underapproximation whose quality is determined by its proximity to the original function:

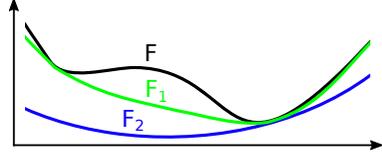

Figure 2. $F_1, F_2$ are convex underapproximations of $F$, where $F_1$ is tighter than $F_2$, *i.e.* $F_1 \geq F_2$. The minimiser of $F_1$ and $F$ coincide, the minimisers of $F_2$ and $F$ do not coincide.

**Definition 1.** *Let $F_1$, $F_2$ be convex underapproximations of $F$. We say that $F_1$ is* at least as tight as $F_2$*, if $F_1(x) \geq F_2(x)$ for all $x$. If in addition $F_1(x) > F_2(x)$ holds for at least one $x$, we say that $F_1$ is* tighter *than $F_2$.*

A systematic way to construct convex underapproximations is to write $F(x) = f_1(x) + f_2(x)$ and use the fact that

$$F^{**}(x) \geq f_1^{**}(x) + f_2^{**}(x), \quad (4)$$

where the *convex conjugate* $f_i^*(p) = \sup_x p^T x - f_i(x)$ of $f_i$ is computed twice to obtain the *biconjugate* $f_i^{**}$.

**Faithful convex underapproximations:** Next, we introduce the notion of *faithful* convex underapproximations.

**Definition 2.** *Let $G : \mathbb{R}^n \to \mathbb{R}$ be a cost function and $C \subset \mathbb{R}^n$ be the feasible set. We call a convex underapproximation $F_{conv}$ of $G + \delta_C$ faithful, if $F_{conv}(x) = G(x)$ holds for all $x \in C$.*

Note that $G$ can be convex or non-convex, and $C$ must not necessarily be convex. Faithful convex underapproximations inherit some obvious but appealing properties:

**Proposition 3.** *If a minimiser $\tilde{x}$ of a faithful convex underapproximation $F_{conv}$ of $G + \delta_C$ meets $\tilde{x} \in C$, it is also a global minimiser of $G + \delta_C$.*

*Proof.* Assume that our statement is false. Then there exists a $x \in C$ with $G(x) < G(\tilde{x})$. However, since $F_{conv}$ is faithful, we have $G(\tilde{x}) = F_{conv}(\tilde{x})$ as well as $G(x) = F_{conv}(x)$. We conclude $F_{conv}(x) < F_{conv}(\tilde{x})$ which contradicts $\tilde{x}$ being a minimiser of $F_{conv}$. ∎

**Proposition 4.** *If $F_{conv}^1$ and $F_{conv}^2$ are two convex underapproximations of $G + \delta_C$, where $F_{conv}^1$ is faithful and $F_{conv}^2$ is not, then $F_{conv}^2$ cannot be tighter than $F_{conv}^1$.*

*Proof.* Since $F_{conv}^2$ is not faithful, there exists a $x \in C$ for which $F_{conv}^2(x) \neq G(x)$. Because $F_{conv}^2$ is still an underapproximation of $G + \delta_C$, it must hold that $F_{conv}^2(x) < G(x) = F_{conv}^1(x)$, which shows that $F_{conv}^2$ cannot be tighter than $F_{conv}^1$. ∎

## 3. Problem Statement

Throughout the rest of the paper we consider the problem of finding a permutation matrix $X \in \mathbb{P}_n$ that minimises quadratic costs, *i.e.*

$$\arg\min_{X \in \mathbb{R}^{n \times n}} f(x) + \delta_\mathbb{P}(x) =: F(x), \quad (5)$$

where $f(x) := x^T W x + c^T x$, and $\delta_\mathbb{P}(x)$ is the indicator function for the set of permutations $\mathbb{P}_n$. Since partial permutations can also be handled in Problem (5) by introducing dummy variables, we assume (full) permutation matrices $X \in \mathbb{P}_n$. Moreover, for the sake of a simpler explanation we assume that $\mathcal{C} = \mathbb{R}^{n \times n}$ in Problem (1).

## 4. A General Convex Relaxation Framework

Before we proceed with details, we give a brief summary how we obtain our convex relaxation: First, we show that there exists an infinite number of functions $\tilde{f}_\Delta$ parametrised by $\Delta$, for which $f(x) = \tilde{f}_\Delta(x)$ holds whenever $\text{vec}^{-1}(x)$ is a permutation matrix. This observation is enormously helpful since it allows us to replace $f$ in (5) with $\tilde{f}_\Delta$ (for a suitably chosen $\Delta$) and then consider a convex relaxation thereof. Subsequently, we analyse the behaviour of $\tilde{f}_\Delta$ over the convex hull of the set $\mathbb{P}$, the set of doubly-stochastic matrices $\mathbb{DS}$, from which we derive how to choose $\Delta$ to obtain a tight convex underapproximation (Sec. 5).

### 4.1. The Energy over $\mathbb{P}$

A key insight for developing tight convex relaxations of (5) is that $f$ is not unique [41, 42, 13, 8]. More precisely, with

$$\tilde{f}(x; D_1, D_2, d) := x^T(W - Z(D_1, D_2, d))x \quad (6)$$
$$+ (c+d)^T x + \langle \mathbf{I}_n, D_1 + D_2 \rangle,$$
$$Z(D_1, D_2, d) := D_1 \otimes \mathbf{I}_n + \mathbf{I}_n \otimes D_2 + \text{diag}(d), \quad (7)$$

where we also write $\tilde{f}_\Delta := \tilde{f}(\cdot; D_1, D_2, d)$ and $Z(\Delta) := Z(D_1, D_2, d)$ with $\Delta := (D_1, D_2, d)$, one finds the following equivalences of energies:

**Proposition 5.** *Let $X \in \mathbb{P}_n$ and $x = \text{vec}(X)$. For any $d \in \mathbb{R}^{n^2}, D_1, D_2 \in \mathbb{R}^{n\times n}$ it holds that $f(x) = \tilde{f}(x; D_1, D_2, d)$.*

*Proof.* With $X \in \mathbb{P}_n$, we have $X^T X = \mathbf{I}_n$ and therefore

$$\langle \mathbf{I}_n, D_1 \rangle = \langle X^T X, D_1 \rangle, \tag{8}$$
$$= \langle X D_1, X \rangle = \langle \text{vec}(X D_1), \text{vec}(X) \rangle, \tag{9}$$
$$= \langle (D_1^T \otimes \mathbf{I}_n) x, x \rangle, \tag{10}$$
$$= x^T (D_1 \otimes \mathbf{I}_n) x. \tag{11}$$

Similarly $\langle \mathbf{I}_n, D_2 \rangle = x^T (\mathbf{I}_n \otimes D_2) x$. Moreover, since $X \in \mathbb{P}_n$, we have that $X_{ij} = X_{ij}^2$ for all $i, j = 1, \ldots, n$. Thus, $x^T \text{diag}(d) x = d^T x$. Combining the above shows that $f(x) - \tilde{f}(x; D_1, D_2, d) = x^T(D_1 \otimes \mathbf{I}_n + \mathbf{I}_n \otimes D_2 + \text{diag}(d))x - d^T x - \langle \mathbf{I}_n, D_1 + D_2 \rangle = 0$. ∎

Since one can make $\tilde{f}_\Delta$ convex by choosing $\Delta$ such that $W - Z(\Delta)$ is positive semi-definite (PSD), in which case $\tilde{f}_\Delta$ and $\tilde{f}_\Delta^{**}$ coincide, the non-convex permutation constraint is the mere reason for the non-convexity of Problem (1).

### 4.2. The Energy over $\mathbb{DS}$

We now present a result of the behaviour of the energy function when the constraint set $\mathbb{P}$ is replaced by its convex hull, the set of doubly-stochastic matrices $\mathbb{DS}$:

**Lemma 6.** *Let $D, D' \in \mathbb{R}^{n \times n}$ be symmetric and let $d, d' \in \mathbb{R}^{n^2}$. Define $\hat{D} = D' - D$. If $d_i \leq d'_i$ for all $i = 1, \ldots, n^2$, as well as $\hat{D}_{ii} - \max_{j \neq i}(\max(\hat{D}_{ij}, 0)) \geq 0$ for all $i = 1, \ldots, n$, then it holds for all $x \in \text{vec}(\mathbb{DS}_n)$*

$$\tilde{f}(x; D, \bullet, \circ) \leq \tilde{f}(x; D', \bullet, \circ), \tag{12}$$
$$\tilde{f}(x; \diamond, D, \circ) \leq \tilde{f}(x; \diamond, D', \circ), \text{ and} \tag{13}$$
$$\tilde{f}(x; \diamond, \bullet, d) \leq \tilde{f}(x; \diamond, \bullet, d'), \tag{14}$$

*where $\diamond \in \mathbb{R}^{n \times n}, \bullet \in \mathbb{R}^{n \times n}$ and $\circ \in \mathbb{R}^{n^2}$.*

*Proof.* We have

$$\tilde{f}(x; D, \bullet, \circ) - \tilde{f}(x; D', \bullet, \circ) \tag{15}$$
$$= x^T(-Z(D, \bullet, \circ) + Z(D', \bullet, \circ))x - \langle \mathbf{I}_n, \hat{D} \rangle \tag{16}$$
$$= x^T(D' \otimes \mathbf{I}_n - D \otimes \mathbf{I}_n)x - \langle \mathbf{I}_n, \hat{D} \rangle \tag{17}$$
$$= x^T((D' - D) \otimes \mathbf{I}_n)x - \langle \mathbf{I}_n, \hat{D} \rangle \tag{18}$$
$$= \langle X^T X, \hat{D} \rangle - \langle \mathbf{I}_n, \hat{D} \rangle \tag{19}$$
$$= \sum_i \left( ((X^T X)_{ii} - 1)\hat{D}_{ii} + \sum_{j \neq i}(X^T X)_{ij} \hat{D}_{ij} \right). \tag{20}$$

We continue by looking at (20) for each $i$ separately:

$$((X^T X)_{ii} - 1)\hat{D}_{ii} + \sum_{j \neq i}(X^T X)_{ij} \hat{D}_{ij} \tag{21}$$
$$\leq ((X^T X)_{ii} - 1)\hat{D}_{ii} + \sum_{j \neq i}(X^T X)_{ij} \max(\hat{D}_{ij}, 0) \tag{22}$$
$$\leq (1 - (X^T X)_{ii})(-\hat{D}_{ii})$$
$$+ \left( \max_{j \neq i}(\max(\hat{D}_{ij}, 0)) \right) \sum_{j \neq i}(X^T X)_{ij} \tag{23}$$
$$= \underbrace{(1 - (X^T X)_{ii})}_{\geq 0} \underbrace{\left( \max_{j \neq i}(\max(\hat{D}_{ij}, 0)) - \hat{D}_{ii} \right)}_{\leq 0 \text{ by assumption}} \tag{24}$$
$$\leq 0. \tag{25}$$

In the step from (23) to (24) we used that if $X$ is doubly-stochastic, then so is $X^T X$. Thus, using (20) it follows that $\tilde{f}(x; D, \bullet, \circ) \leq \tilde{f}(x; D', \bullet, \circ)$. The case in (13) is analogous. Note that tighter (but more complicated criteria) can be derived by additionally considering the sum over $i$ and using that $(X^T X)_{i,i} \geq \frac{1}{n}$. We skipped this analysis for the sake of simplicity.

Moreover,

$$\tilde{f}(x; \diamond, \bullet, d) - \tilde{f}(x; \diamond, \bullet, d') \tag{26}$$
$$= x^T(-Z(\diamond, \bullet, d) + Z(\diamond, \bullet, d'))x + (d - d')^T x \tag{27}$$
$$= x^T(\text{diag}(d') - \text{diag}(d))x + (d - d')^T x \tag{28}$$
$$= x^T(\text{diag}(d' - d))x + (d - d')^T x \tag{29}$$
$$= \sum_{i=1}^{n^2}(d' - d)_i x_i^2 + (d - d')_i x_i \tag{30}$$
$$= \sum_{i=1}^{n^2}(d' - d)_i x_i^2 - (d' - d)_i x_i \tag{31}$$
$$= \sum_{i=1}^{n^2}(d' - d)_i (x_i^2 - x_i) \leq 0. \tag{32}$$

The last inequality follows from the assumption $d'_i - d_i \geq 0$ and $x_i^2 - x_i \leq 0$ (using $x \in \text{vec}(\mathbb{DS}_n)$). ∎

For the special case of diagonal matrices $D_1$ and $D_2$, Lemma 6 indicates that larger entries as well as larger elements of $d$ lead to tighter relaxations.

## 5. Tight Convex Relaxations

The previous sections suggest the following strategy:
1. To obtain *tight* relaxations we want to write $F = f_1 + \delta_\mathbb{P}$, and approximate $F^{**}$ by $f_1^{**} + \delta_{\mathbb{DS}_n}$ in such a way that $f_1^{**}$ is as large as possible.

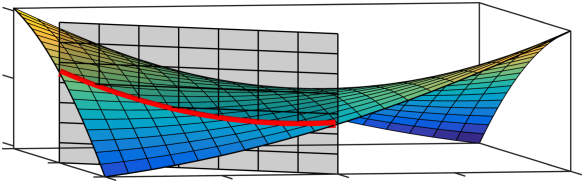

Figure 3. Illustration of subspace convexity. The non-convex quadratic energy function (coloured surface) becomes convex (red line) when restricted to a subspace (grey plane).

2. To obtain a *faithful* relaxation, $f_1$ needs to be convex. Unfaithful relaxations cannot be tighter than faithful ones.

Based on the above demands, a naive choice would therefore be $f_1(\cdot; \Delta) := \tilde{f}(\cdot; \Delta)$ for $\Delta$ such that $W - Z(\Delta)$ is positive semi-definite, and which is optimal in the sense of Lemma 6. However, one can obtain tighter relaxations by additionally restricting $f_1$ to an affine subspace [8, 17].

### 5.1. Subspace Convexity

The sets $\mathbb{P}_n$ and $\mathbb{DS}_n$ are subsets of the affine subspace

$$\mathcal{A} = \{x \in \mathbb{R}^{n^2} \ : \ Ax = \mathbf{1}_{2n}\} \text{ for } A = \begin{bmatrix} \mathbf{I}_n \otimes \mathbf{1}_n^T \\ \mathbf{1}_n^T \otimes \mathbf{I}_n \end{bmatrix}, \quad (33)$$

*i.e.* the set of $x = \text{vec}^{-1}(X)$ for matrices $X$ whose rows and columns sum to one. Therefore, we propose to consider the (redundant) splitting

$$F(x) = \underbrace{\tilde{f}(x; \Delta) + \delta_{\mathcal{A}}(x)}_{f_1} + \delta_{\mathbb{P}}(x),$$

where $\delta_{\mathcal{A}}$ is the indicator function of the set $\mathcal{A}$. Note that $f_1$ is convex whenever $W - Z(\Delta)$ is PSD on the subspace $\mathcal{A}$, *i.e.* whenever $F^T(W - Z(\Delta))F \succeq 0$, where $F \in \mathbb{R}^{n^2 \times (n^2 - 2n + 1)}$ denotes a matrix that spans $\ker(A)$. This is a strictly weaker condition than the PSDness of $W - Z(\Delta)$, as we have illustrated in Fig. 3. Note that if $\mathcal{C}$ is (a subset of) a subspace, one can additionally take this into account to obtain an even tighter relaxation.

### 5.2. The Proposed DS* Relaxation

We introduce our convex relaxation of Problem (5) as

$$(\text{DS*}) \qquad \arg\min_{X \in \mathbb{DS}_n} \tilde{f}(x; \Delta), \qquad (34)$$

for a suitable choice $\Delta$ such that $\tilde{f}(x; \Delta)$ is convex.

**Special cases:** The *DS++* approach [17] is obtained as special case of DS* by choosing $\Delta_{\text{DS++}} = (\lambda_{\min}^\star \mathbf{I}_n, \mathbf{0}, \mathbf{0})$, where $\lambda_{\min}^\star$ is the smallest eigenvalue of $F^T W F$ (with $F^T F = \mathbf{I}$), such that the convexity is only enforced on the subspace $\mathcal{A}$. The *DS+* relaxation [20] aims for convexity on the entire $\mathbb{R}^{n^2}$ and is obtained using $\Delta_{\text{DS+}} = (\lambda_{\min} \mathbf{I}_n, \mathbf{0}, \mathbf{0})$, for $\lambda_{\min}$ being the smallest eigenvalue of $W$.

**Proposed relaxation:** Among all possible relaxations of the form (34), we propose to choose the relaxation that maximises the energy in the *centroid* of the simplex of doubly-stochastic matrices, i.e. for $x = \frac{1}{n}\mathbf{1}_{n^2}$, as this is the most undesirable point when truly seeking for permutation matrices. Based on Lemma 6, we add constraints that ensure that the solution is at least as tight as *DS++* while maximising the value at $x = \frac{1}{n}\mathbf{1}_{n^2}$. We point out that there are less restrictive constraints, which, however, are more difficult to interpret than the ones we have chosen to present below.

**Proposition 7.** *The minimiser $\tilde{\Delta}$ among all $\Delta = (D_1, D_2, d)$ with symmetric $D_1$ and $D_2$ of*

$$\min_{\Delta} \quad -\text{tr}(Z(\Delta)) + \frac{1}{n}\sum_{i,j}(Z(\Delta))_{ij} \qquad (35)$$

$$\text{s.t.} \quad F^T(W - Z(\Delta))F \succeq 0,$$
$$0 \geq -(D_1)_{ii} + \max_{j \neq i} \max((D_1)_{ij}, 0) \quad \forall\, i,$$
$$0 \geq -(D_2)_{ii} + \max_{j \neq i} \max((D_2)_{ij}, 0) \quad \forall\, i,$$
$$d_i \geq \lambda_{\min}^\star \quad \forall\, i,$$

*yields a relaxation that is at least as tight as* DS++. *If $Z(\tilde{\Delta}) \neq \lambda_{\min}^\star \mathbf{I}_{n^2}$, the above is tighter than* DS++.

*Proof.* First of all, we observe that, in addition to $\Delta_{\text{DS++}} = (\lambda_{\min}^\star \mathbf{I}_n, \mathbf{0}, \mathbf{0})$, the DS++ relaxation is also obtained by the choice $\Delta'_{\text{DS++}} = (\mathbf{0}, \mathbf{0}, \lambda_{\min}^\star \mathbf{1}_{n^2})$. The constraints in (35) are feasible as they are satisfied for the *DS++* choice $\Delta'_{\text{DS++}} = (\mathbf{0}, \mathbf{0}, \lambda_{\min}^\star \mathbf{1}_{n^2})$. Thus, a minimiser exists.

Writing $\tilde{\Delta} = (\tilde{D}_1, \tilde{D}_2, \tilde{d})$, the convex constraints immediately yield that

$$\tilde{f}(x; \lambda_{\min}^\star \mathbf{I}_n, \mathbf{0}, \mathbf{0}) = \tilde{f}(x; \mathbf{0}, \mathbf{0}, \lambda_{\min}^\star \mathbf{1}_{n^2})$$
$$\leq \tilde{f}(x; \mathbf{0}, \mathbf{0}, \tilde{d})$$
$$\leq \tilde{f}(x; \tilde{D}_1, \mathbf{0}, \tilde{d})$$
$$\leq \tilde{f}(x; \tilde{D}_1, \tilde{D}_2, \tilde{d}) = f(x; \tilde{\Delta})$$

holds for all $x \in \text{vec}(\mathbb{DS}_n)$ based on Lemma 6.

Finally, one can compare $\tilde{f}(x; \tilde{\Delta})$ with $\tilde{f}(x; \Delta'_{\text{DS++}})$ at $x = \frac{1}{n}\mathbf{1}_{n^2}$ to see that the *DS++* relaxation is strictly below the relaxation given by (35) if *DS++* does not happen to yield a solution to (35) already. ∎

### 5.3. Efficient Approximations

Although the proposed approach in (35) has the three advantages that (i) the solution is optimal in some sense, (ii) it can be computed using convex optimisation techniques, and (iii) we only have $3n^2$ unknowns, the semi-definite constraint $F^T(W - Z(\Delta))F \succeq 0$ involves a large matrix of size $(n^2 - 2n + 1) \times (n^2 - 2n + 1)$. To enforce the latter using first-order methods, one needs to iteratively project onto the PSD cone of such matrices, which has similar complexity as the lifting-based relaxation approaches [59, 27]. The key question for practical applications therefore becomes how to approximate (35) such that the resulting algorithm is scalable.

To do so we consider the special case of $D_1$ and $D_2$ being diagonal matrices and $d=\mathbf{0}_{n^2}$, leaving us with $2n$-many unknowns. Note that in this case we still have a relaxation that is at least as tight as *DS++* (cf. the paragraph about the special case in Sec. 5.2). While the constraints in (35) have to be modified for such a choice, we keep in mind that their purpose is to prevent individual entries from becoming too small, and focus on reducing the computational costs of handling the PSD constraint. We follow some ideas of [52, 53], where such a constraint is replaced by the penalty

$$\tilde{h}(Y) = \|Y_-\|_F^2 = \sum_i \min(\lambda_i(Y), 0)^2, \quad (36)$$

for $\lambda(Y)$ denoting the spectrum of $Y$. Since the gradient evaluation of $\tilde{h}$ requires the computation of all negative eigenvalues of $Y$, we propose to *only penalise the smallest negative eigenvalue*. We define

$$h(Y) = \frac{1}{2}\min(\lambda_{\min}(Y), 0)^2, \quad (37)$$
$$T(d_1, d_2) = F^T(W - \operatorname{diag}(d_1) \otimes \mathbf{I}_n - \mathbf{I}_n \otimes \operatorname{diag}(d_2))F,$$

and introduce $h(T(d_1, d_2))$ into our objective, where $d_1$ and $d_2$ denote the diagonals of $D_1$ and $D_2$. Interestingly, $h(T(d_1, d_2))$ is differentiable if the smallest eigenvalue of $T(d_1, d_2)$ has multiplicity one [32]. Moreover, gradients can be computed efficiently due to two reasons: (i) A gradient of $h$ merely requires the smallest eigenvalue/eigenvector pair which can be determined via an inverse power method. (ii) The special structure of the adjoint of the affine operator $T$ allows to efficiently compute the inner derivative:

**Lemma 8.** *Let $T(d_1, d_2)$ have a smallest eigenvalue $\lambda_{\min}$ of multiplicity 1, and let $u_{\min}$ be a corresponding eigenvector with $\|u_{\min}\| = 1$. Then*

$$(p_1)_j = -\min(\lambda_{\min}, 0) \sum_i ((\operatorname{vec}^{-1}(Fu_{\min}))_{i,j})^2,$$
$$(p_2)_i = -\min(\lambda_{\min}, 0) \sum_j ((\operatorname{vec}^{-1}(Fu_{\min}))_{i,j})^2$$

*meet $p_1 = \nabla_{d_1}(h \circ T)(d_1, d_2)$, $p_2 = \nabla_{d_2}(h \circ T)(d_1, d_2)$.*

*Proof.* The proof is based on the fact that $(h \circ T)$ is a composition of four functions:

$$(h \circ T)(d_1, d_2) = \frac{1}{2}\min(g(\lambda(T(d_1, d_2))), 0)^2,$$

*i.e.* the affine function $T$, a function $Y \mapsto \lambda(Y)$ determining the eigenvalues of $Y$, a function $g(v) = \min(v)$ selecting the minimal element of a vector, and the function $x \mapsto \frac{1}{2}\min(x, 0)^2$. The latter is continuously differentiable with derivative $\min(x, 0)$.

Compositions of the form $g(\lambda(Y))$ have been studied in detail in [32], and according to [32, p. 585, Example of Cox and Overton] it holds that

$$\partial(g \circ \lambda)(Y) = \operatorname{conv}\{uu^T \ : \ Yu = \lambda_{\min}(Y)u, \|u\| = 1\}. \quad (38)$$

Note that $(g \circ \lambda)$ becomes differentiable if the smallest eigenvalue of $Y$ has multiplicity one, such that the corresponding eigenspace is of dimension 1 and the above set $\partial(g \circ \lambda)(Y)$ reduces to a singleton – also see [32, Theorem 2.1].

Thus, by the chain rule

$$\min(\lambda_{\min}(Y), 0) \, u_{\min} u_{\min}^T$$

is a gradient of $h$ at $Y$ if $\lambda_{\min}(Y)$ has multiplicity 1.

Left to consider is the inner derivative coming from the affine map $T$. Let us consider the linear operator

$$\tilde{T}(d_1) = -F^T(\operatorname{diag}(d_1) \otimes \mathbf{I}_n)F$$

as the part of $T$ that has a relevant inner derivative with respect to $d_1$. The gradient of a linear operator $\tilde{T}$ is nothing but its adjoint operator $\tilde{T}^*$, *i.e.* the operator for which

$$\langle \tilde{T}(d), A \rangle = \langle d, \tilde{T}^*(A) \rangle$$

holds for all $d$ and all $A$. (In this case we could explicitly prove this by vectorizing the entire problem, but the relation holds in much more generality as the definition of general (Gateaux) gradients utilises the Riesz representation theorem, see e.g. [3, p. 40, Remark 2.55]). Since the adjoint of $T_1 \circ T_2$ is $T_2^* \circ T_1^*$, we can consider the operations separately in a reverse order. The last thing $\tilde{T}$ does is the multiplication with $F^T$ from the left and with $F$ from the right, which means that the first thing the adjoint $\tilde{T}^*$ does is the multiplication with $F$ from the left and with $F^T$ from the right.

The operator $\operatorname{diag}(d_1) \otimes \mathbf{I}_n$ repeats the entries of $d_1$ $n$ times, and writes the result on the diagonal of an $n^2 \times n^2$ diagonal matrix. The adjoint of writing a vector of length $n^2$ on the diagonal of an $n^2 \times n^2$ diagonal matrix, is the extraction of the diagonal of such a matrix. Finally, the adjoint of the repeat operation is the summation over the components of those indices at which values were repeated. As an illustrative example, note that

$$\underbrace{A = \begin{pmatrix} 1 & 0 \\ 1 & 0 \\ 0 & 1 \\ 0 & 1 \end{pmatrix}}_{\text{repeat each component}} \Rightarrow \underbrace{A^* = \begin{pmatrix} 1 & 1 & 0 & 0 \\ 0 & 0 & 1 & 1 \end{pmatrix}}_{\text{sum over repeated components}}.$$

If $\tilde{T}^*$ is applied to an element $Y = uu^T \in \mathbb{R}^{n^2 \times n^2}$ the first steps are left multiplication with $F$ and right multiplication with $F^T$, leading to $(Fu)(Fu)^T$. The extraction of the diagonal of the resulting matrix yields a vector of length $n^2$ with entries $(Fu)_k^2$. By taking sums over $n$ consecutive entries, and multiplying with the remaining inner derivatives $(-1)$ and $\min(\lambda_{\min}, 0)$ we arrive at the formula for $\nabla_{d_1}(h \circ T)$ as stated by Lemma 8. Determining the formula for $\nabla_{d_2}(h \circ T)$ follows exactly the same computation with a different final summation as the operator $\mathbf{I}_n \otimes \text{diag}(d_2)$ repeats the entries in a different order. ■

Considering the great success of subgradient descent type of methods in computer vision, e.g. in the field of deep learning, we consider such a method for our problem, too: We simply keep using the formulas in Lemma 8 even in the case where the multiplicity of the smallest eigenvalue is larger than 1, and just select any $u_{\min}$.

Our general strategy is to take a formulation like (35), replace the semi-definite constraint by the penalty $h(T(d_1, d_2))$ in (37), run a few iterations of a subgradient descent type algorithm on the resulting energy to obtain $(d_1, d_2)$, and minimise $\tilde{f}(x, \text{diag}(d_1), \text{diag}(d_2), \mathbf{0}_{n^2})$ over $\mathbb{DS}_n$ to obtain the solution $X$ to our convex relaxation. Before we present our algorithm to determine $(d_1, d_2)$, we address the projection of a solution $X \notin \mathbb{P}_n$ onto the set $\mathbb{P}_n$, which imposes some additional demands upon $(d_1, d_2)$.

## 6. Projection onto $\mathbb{P}_n$

Since the objective of the approximation in Sec. 5.3 has the form of Prop. 5, it is a *faithful* convex underapproximation and thus Prop. 3 is applicable. Hence, if the obtained solution $X$ of our relaxation is in $\mathbb{P}_n$, then $X$ is a global solution to Problem (5). However, if $X \notin \mathbb{P}_n$, a strategy for projecting $X$ onto $\mathbb{P}_n$ is necessary. Whilst the $\ell_2$-projection is straightforward, it makes the (over-simplified) assumption that the sought global solution $\hat{X}$ is the permutation matrix that is closest to $X$ in the Euclidean sense.

Instead, we pursue a convex-to-concave path-following, as discussed in Sec. 1.1. To ensure that the final problem we solve is concave, we not only impose PSDness upon $T(d_1, d_2)$, but also NSDness upon $T(-d_1, -d_2)$, so that

$$T_\alpha = (1 - \alpha)T(d_1, d_2) + \alpha T(-d_1, -d_2) \quad (39)$$

represents a path from a PSD ($\alpha=0$) to a NSD ($\alpha=1$) matrix [9]. The semi-definite constraints automatically constrain the matrix $Z(d_1, d_2, \mathbf{0}_{n^2})$ to be non-positive. Moreover, as discussed in Section 5.3, one wants to maximise the sum of the diagonal elements while simultaneously restricting them from becoming too small. In our numerical experiments we found that a quadratic penalty on the diagonal elements is an easy-to-compute way of achieving this. By again replacing the hard semi-definite constraints by their respective cheap soft-constraints we end up with a problem of the form

$$\min_{d_1, d_2} \frac{\eta}{2}(\|d_1\|^2 + \|d_2\|^2) + (1 - \beta)h(T(d_1, d_2)) \quad (40)$$
$$+ \beta h(-T(-d_1, -d_2)).$$

To find $(d_1, d_2)$ we optimise (40) using a proximal subgradient type method with a subgradient type step on the $h$ penalties (as discussed in Sec. 5.3), followed by a proximal step with respect to the quadratic regularisation, see Alg. 1.

Since the PSD and NSD constraints are modelled as penalties, there is no guarantee that the resulting $(d_1, d_2)$ lead to PSD and NSD matrices $T(d_1, d_2)$ and $T(-d_1, -d_2)$, respectively. To compensate for this, one can shift the diagonal of the matrices $T(d_1, d_2)$ and $T(-d_1, -d_2)$ by their smallest/largest eigenvalues, similar to DS++ [17].

Defining $\Delta_\alpha := (d_1 - 2\alpha d_1, d_2 - 2\alpha d_2, \mathbf{0})$, our PF procedure starts with minimising $\tilde{f}_{\Delta_\alpha}$ for $\alpha=0$ over $\mathbb{DS}_n$, which corresponds to the convex Problem (34) that can be solved to global optimality. Then, we gradually increase $\alpha$ and (lo-

**Input:** $W, F, n, n_{\text{iter}}=10, \tau > 0, \eta > 0, \beta \in [0, 1]$
**Output:** $d_1, d_2$
**Initialise:** $d_1 = \mathbf{0}, d_2 = \mathbf{0}$
1 **foreach** $i = 1, \ldots, n_{iter}$ **do**
2 $\quad$ compute $T_0$ and $T_1$ using (39)
3 $\quad$ $[u_{\min}, \lambda_{\min}] = \text{eigs}(T_0, 1, \text{'sa'})$
4 $\quad$ $[u_{\max}, \lambda_{\max}] = \text{eigs}(T_1, 1, \text{'la'})$
5 $\quad$ $V^+ = \text{reshape}((Fu_{\min}) \odot (Fu_{\min}), n, n)$
6 $\quad$ $V^- = \text{reshape}((Fu_{\max}) \odot (Fu_{\max}), n, n)$
$\quad$ // gradient step
7 $\quad$ $d_1 = d_1 + (1-\beta)\tau\lambda_{\min}(V^+)^T\mathbf{1}_n - \beta\tau\lambda_{\max}(V^-)^T\mathbf{1}_n$
8 $\quad$ $d_2 = d_2 + (1-\beta)\tau\lambda_{\min}V^+\mathbf{1}_n - \beta\tau\lambda_{\max}V^-\mathbf{1}_n$
$\quad$ // proximal step of squared $\ell_2$-norm
9 $\quad$ $d_1 = \frac{1}{1+\tau\eta}d_1; \quad d_2 = \frac{1}{1+\tau\eta}d_2$

**Algorithm 1:** Proximal subgradient descent type of algorithm to find $(d_1, d_2)$. The notation eigs(), reshape() is borrowed from MATLAB and $\odot$ denotes the Hadamard product. The parameters $\tau$, $\eta$ and $\beta$ are the step size, the regularisation weight, and the relative importance of $T(d_1, d_2)$ being PSD and $T(-d_1, -d_2)$ being NSD.

cally) minimise $\tilde{f}_{\Delta_\alpha}$ over $\mathbb{DS}_n$ using the previous solution as initialisation. Once $\alpha=1$, the minimisation of the concave function $\tilde{f}_{\Delta_\alpha}$ over $\mathbb{DS}_n$ results in a solution $X \in \mathbb{P}_n$. We use the Frank-Wolfe (FW) method [21] for the minimisation of all $\tilde{f}_{\Delta_\alpha}$ over $\mathbb{DS}_n$, where the linear programming subproblems are solved via the Auction algorithm [7, 6]. Our DS* projection based on PF is shown in Alg. 2.

## 7. Complexity Analysis

Computing $T_0$ and $T_1$ in line 2 in Alg. 1 involves two matrix products with matrices of size $\mathcal{O}(n^2 \times n^2)$. However, the matrix $F \in \mathbb{R}^{n^2 \times (n^2 - 2n + 1)}$ spanning the null space of $A$ in (33) can be constructed as sparse matrix:

**Input:** $W, \Delta_\alpha$ (obtained from Alg. 1)
**Output:** $X \in \mathbb{P}_n$
1 **for** $\alpha = 0, \ldots, 1$ **do**
2 $\quad X = \text{frankWolfe}(\tilde{f}_{\Delta_\alpha}, X)$

**Algorithm 2:** DS* projection. frankWolfe$(f, X_0)$ finds a local minimiser of $f$ over $\mathbb{DS}_n$ with $X_0$ as initialisation.

**Lemma 9.** *Let $x^i \in \mathbb{R}^n$ be defined as*

$$(x^i)_j := \begin{cases} 1 & \text{if } j = i \\ -1 & \text{if } j = i+1 \\ 0 & \text{otherwise} \end{cases}, \quad \text{and let} \quad z^{i,j} := x^i \otimes x^j.$$

*With $F = [z^{1,1}, z^{1,2}, \ldots, z^{n-1,n-1}] \in \mathbb{R}^{n^2 \times (n-1)^2}$, we have that $\text{im}(F) = \ker(A)$ for $A = \begin{bmatrix} \mathbf{I}_n \otimes \mathbf{1}_n^T \\ \mathbf{1}_n^T \otimes \mathbf{I}_n \end{bmatrix}$.*

*Proof.* The linear independence of all $x^1, \ldots x^{n-1}$ implies the linear independence of all $z^{i,j} = x^i \otimes x^j$ for $i, j \in [n-1]$, from which we see that $\dim(\text{im}(F)) = \text{rank}(F) = (n-1)^2 = n^2 - 2n + 1 = \dim(\ker(A))$.

We proceed by showing that $\text{im}(F) \subseteq \ker(A)$. Let $z \in \text{im}(F)$, so $z = \sum_{i,j=1}^{n-1} a_{ij} z^{i,j}$ for some coefficients $\{a_{ij} \in \mathbb{R}\}$. By construction of the $z^{i,j}$, for $i, j \in [n-1]$ we have

$$(\mathbf{I}_n \otimes \mathbf{1}_n^T) z^{i,j} = \mathbf{0}_n \text{ and } (\mathbf{1}_n^T \otimes \mathbf{I}_n) z^{i,j} = \mathbf{0}_n, \quad (41)$$

which implies that

$$(\mathbf{I}_n \otimes \mathbf{1}_n^T) a_{ij} z^{i,j} = \mathbf{0}_n \text{ and } (\mathbf{1}_n^T \otimes \mathbf{I}_n) a_{ij} z^{i,j} = \mathbf{0}_n. \quad (42)$$

Thus

$$(\mathbf{I}_n \otimes \mathbf{1}_n^T) z = \mathbf{0}_n \text{ and } (\mathbf{1}_n^T \otimes \mathbf{I}_n) z = \mathbf{0}_n, \quad (43)$$

from which we can see that $z \in \ker(A)$. Combining $\dim(\text{im}(F)) = \dim(\ker(A))$ and $\text{im}(F) \subseteq \ker(A)$ shows that $\text{im}(F) = \ker(A)$. ∎

As $F$ is a sparse matrix with $\mathcal{O}(n^2)$ non-zero elements, computing $T_0$ and $T_1$ in line 2 has complexity $\mathcal{O}(n^4)$. The eigenvalue computations in lines 3 and 4 have complexity $\mathcal{O}(n^4)$ when running an iterative solver for a fixed number of iterations. Thus, Alg. 1 has complexity $\mathcal{O}(n^4)$, which is $\mathcal{O}(n^2)$ smaller than the lifted SDP relaxations [59, 27].

## 8. Experiments

To evaluate DS*, we first compare bounds obtained by (lifting-free) convex relaxation methods on synthetic matching problems. Then, we consider image arrangement using three different datasets. Eventually, we incorporate our approach into a convex multi-graph matching method which we evaluate on synthetic and real data. Unless stated otherwise we used 10 steps for PF for DS+, DS++ and DS*. Moreover, for DS* we used $\tau = 4, \eta = 0.1$ and $\beta = 0.2$.

### 8.1. Synthetic Data

Here, we compare DS+, DS++ and DS* on random instances of Problem (5). To do so, for each $n \in \{16, 20, \ldots, 40\}$ we draw 200 symmetric matrices $W \in \mathbb{R}^{n \times n}$ with uniformly distributed elements in $]-1, 1[$. In total, we solve $7 \cdot 200$ optimisation problems for each method. The lower bounds are given by the objective value of the respective relaxation method, and the upper bounds are obtained via projection using PF. To allow a comparison across individual problem instances, we normalise the objective values: After solving one problem instance with the three methods, we scale the three objective values such that the largest lower bound is $-1$. Upper bounds are normalised analogously (independent of the lower bounds).

**Results:** The mean and the standard deviation of the so-obtained lower and upper bounds are shown in Fig. 4. It can be seen that our proposed approach results in better lower bounds compared to DS+ and DS++, as expected due to our theoretical results. In general, better lower bounds are no guarantee for better upper bounds, as can be seen when comparing the upper bounds of DS+ and DS++. However, our method is also able to achieve the best upper bounds.

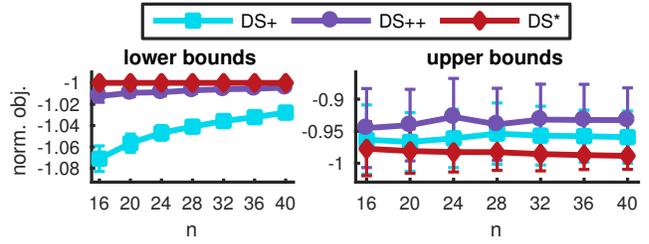

Figure 4. Comparison of *DS+*, *DS++* and *DS\**. Left: lower bounds (higher is better). Right: upper bounds obtained by projection, cf. Sec. 6 (lower is better).

### 8.2. Image Arrangement

In this experiment we consider the arrangement of a collection of images on a predefined grid such that "similar" images are close to each other on the grid (see Fig. 5). In [22] this is tackled by minimising the energy

$$E(X) = \delta_\mathbb{P}(X) + \min_c \sum_{ijkl} X_{ij} X_{kl} |c \cdot d_{ik} - d'_{jl}|, \quad (44)$$

where the scalar factor $c$ is used as normalisation between the pairwise image "distances" $d \in \mathbb{R}^{n \times n}$ and the pairwise grid position distances $d' \in \mathbb{R}^{n \times n}$. The distances $d$ ($d'$) are computed as the $\ell_2$-norm of the differences between pairs of image features (or grid positions). To employ image arrangement with DS++ and DS*, we fix $c$ in (44) such that $d$ and $d'$ have the same mean [17].

**Setup:** We compare DS* with *isomatch* (using random swaps) [22] and DS++ [17] on various datasets (random colours, face images [49], and COCO [33]), where we used the RGB colour, MoFA facial expression parameters [49], and the average hue-saturation vector of each image as features, respectively. For the *random colours* experiment, in each run we uniformly sample random RGB values and then arrange the individual colours on the grid (i.e. in the "image arrangement" terminology we arrange images that comprise a single pixel). For the *face* and *COCO* experiments, in each run we randomly select images that are arranged on a grid.

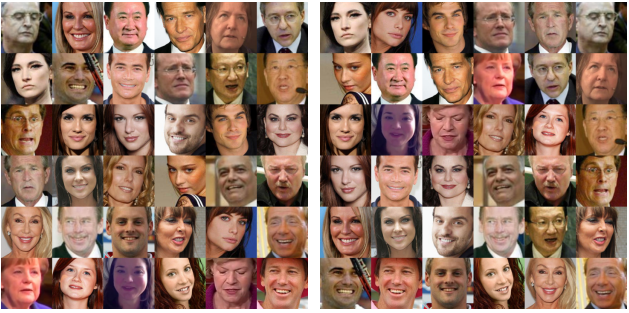

Figure 5. Left: 36 face images are randomly arranged on a $6 \times 6$ grid. Right: The face images are arranged with DS* according to facial expression [49]. (Best viewed on screen when zoomed in)

Table 1. Quantitative image arrangement results. We use $^\star$ to indicate sparse matchings of $\tilde{n}$ elements that are subsequently extrapolated using PMF. Cases that are not indicated by $^\star$ mean that we used all images for directly computing the matching.

| data/<br>feat. | grid<br>size | $\tilde{n}$ | initial | isomatch | DS++ | DS* |
|---|---|---|---|---|---|---|
| rnd. colours RGB | $8^2$ | all | 0.466 | 0.227 | – | 0.211 | **0.196** |
| | $16^2$ | 50 | 0.475 | 0.288 | $0.249^\star$ | $0.242^\star$ | $\mathbf{0.236}^\star$ |
| | $32^2$ | 75 | 0.476 | – | $0.246^\star$ | $0.245^\star$ | $\mathbf{0.235}^\star$ |
| | $64^2$ | 75 | 0.478 | – | $0.259^\star$ | $0.249^\star$ | $\mathbf{0.244}^\star$ |
| face facial exp. | $8^2$ | all | 0.519 | 0.285 | – | 0.279 | **0.266** |
| | $16^2$ | 50 | 0.530 | 0.324 | $0.307^\star$ | $0.305^\star$ | $\mathbf{0.300}^\star$ |
| | $32^2$ | 75 | 0.531 | – | $0.352^\star$ | $0.336^\star$ | $\mathbf{0.324}^\star$ |
| | $42^2$ | 75 | 0.533 | – | $0.379^\star$ | $0.360^\star$ | $\mathbf{0.344}^\star$ |
| COCO avg. hue/sat. | $8^2$ | all | 0.541 | 0.273 | – | 0.256 | **0.244** |
| | $16^2$ | 50 | 0.548 | 0.315 | $0.306^\star$ | $\mathbf{0.238}^\star$ | $\mathbf{0.238}^\star$ |
| | $32^2$ | 75 | 0.552 | – | $0.330^\star$ | $\mathbf{0.274}^\star$ | $\mathbf{0.274}^\star$ |
| | $40^2$ | 75 | 0.549 | – | $0.338^\star$ | $\mathbf{0.291}^\star$ | $0.292^\star$ |

Table 2. Average runtimes for the random colours experiments.
$^\dagger$ Due to the slow processing in the *isomatch* $32^2$ and $64^2$ settings we have only run these settings once to estimate the runtime.

| grid size | $\tilde{n}$ | isomatch | DS++ | DS* |
|---|---|---|---|---|
| $8^2$ | all | 1.29s | – | 17.84s | 43.95s |
| $16^2$ | 50 | 16.36s | $1.11s^\star$ | $5.85s^\star$ | $19.56s^\star$ |
| $32^2$ | 75 | $166.44s^\dagger$ | $2.63s^\star$ | $10.14s^\star$ | $29.53s^\star$ |
| $64^2$ | 75 | $3845.99s^\dagger$ | $31.44s^\star$ | $32.68s^\star$ | $52.83s^\star$ |

**Large-scale arrangement:** In addition to the medium-scale problems of arranging a few hundred images [22, 17], we also investigate the more challenging case of large matching problems that arrange thousands of images. To this end, we first compute a matching between a subset of images and grid cells, which is then extrapolated to obtain a full matching. To be more specific, we use farthest point sampling to obtain the subset of $\tilde{n}$ images as well as $\tilde{n}$ grid cells. Then, using the respective method we solve the (small) arrangement problem between those selected $\tilde{n}$ images and grid cells to obtain a partial matching. Eventually, in order to retrieve a full matching we apply the *Product Manifold Filter* (PMF) method [51] with the computed partial matching as initialisation (we use 5 PMF iterations, and the kernel bandwidth is set to the standard deviation of the elements in $d$ and $d'$, respectively; see [51] for details).

**Results:** Qualitative results for arranging face images according to facial expressions are shown in Fig. 5. Table 1 shows quantitative results for all datasets, where we show the objective value of (44) averaged over 100 runs. It can be seen that the DS* method is on par with DS++ for COCO images, and outperforms the other approaches for random colours and faces. The large-scale image arrangement results show that DS* is a powerful way for initialising PMF. Runtimes are shown in Table 2. To obtain a fair comparison, in the *isomatch* implementation of [22] we replaced the Hungarian [35] LAP solver by the more efficient Auction algorithm [7] as implemented in [6].

### 8.3. Multi-graph Matching

The aim of multi-graph matching (MGM) is to obtain a matching between $k{>}2$ graphs. One way of formulating MGM is to consider all pairwise matchings $\mathbf{X} := [X_{ij}]_{i,j\in[k]} \in (\mathbb{P}_n)^{k\times k}$ and ensure that they are transitively consistent, *i.e.* for all $i,j,\ell \in [k]$ it holds that $X_{ij}X_{j\ell}{=}X_{i\ell}$. With $f^{ij}(x_{ij}) := x_{ij}^T W_{ij} x_{ij}$ being the matching costs between graphs $i$ and $j$, the MGM problem reads

$$\min_{\mathbf{X}\in(\mathbb{P}_n)^{k\times k}} \sum_{i,j\in[k]} f^{ij}(x_{ij}) \quad (45)$$

$$\text{s.t.} \quad X_{ij}X_{j\ell} = X_{i\ell} \ \forall\, i,j,\ell \in [k]\,.$$

We propose to use a convex relaxation of Problem (45):

$$\min_{\mathbf{X}\in(\mathbb{DS}_n)^{k\times k}} \sum_{i,j\in[k]} \tilde{f}^{ij}(x_{ij};\Delta^{ij}) \quad (46)$$

$$\text{s.t.} \quad \mathbf{X} \succeq 0,\ X_{ii} = \mathbf{I}_n \ \forall\, i\,,$$

for $\{\Delta^{ij} : i,j{\in}[k]\}$, where each $\Delta^{ij}$ is determined as in Sec. 5 and 6. The constraints $\mathbf{X} \succeq 0$, $X_{ii} = \mathbf{I}_n$ are the relaxation of the transitive consistency constraint (see [37, 25, 5, 27] for details). While Problem (46) is convex, its objective is *nonlinear* and thus standard SDP solvers are not directly applicable. Instead, we introduce the constraint $\mathbf{X} \succeq 0$ in Problem (46) into the objective function using the soft-penalty $\sigma\|\mathbf{X}_-\|_F^2$, see (36), which we then minimise using the FW method. For solving Problem (46), we conduct PF over the $\tilde{f}^{ij}$, as described in Sec. 6. Since the transitive consistency constraint is relaxed, the resulting solution is not necessarily transitively consistent, which we tackle using *permutation synchronisation* [37].

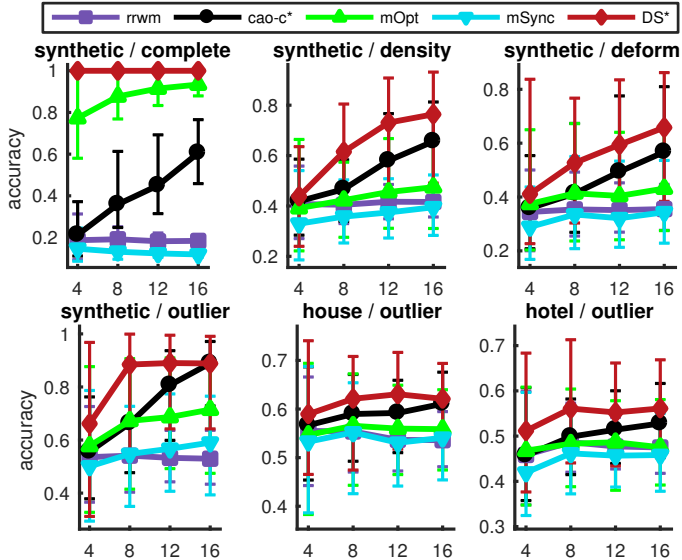

Figure 6. Comparison of accuracies for MGM methods. The vertical lines indicate the root mean square deviation of all the accuracy values above the mean and below the mean, respectively. Each plot shows a different pair of dataset and configuration. The number of graphs $k$ varies along the horizontal axis.

**Setup:** We compare our MGM approach to RRWM [11], composition-based affinity optimisation (CAO) [55], MatchOpt (mOpt) [56], and permutation synchronisation (mSync) [37]. We consider three datasets, synthetic problems and MGM problems using the CMU *house* and *hotel* sequence. For the evaluation we follow the protocol implemented by the authors of [55], where further details are described. We set $\sigma{=}k^{-1}16{,}000$ and we use 30 PF steps.

**Results:** Fig. 6 shows that DS* considerably outperforms the other methods We argue that DS* has superior performance because we simultaneously consider transitive consistency and the matching costs during optimisation. Thus, our approach is better able to leverage the available information. Whilst a related MGM approach has been presented in [27], the authors consider a lifting of the pairwise matching matrices, which is only applicable to very small-scale problems due to the $\mathcal{O}(n^4k^2)$ variables (cf. Fig. 1), in contrast to our approach with only $\mathcal{O}(n^2k^2)$ variables.

## 9. Discussion and Future Work

We have found that running Alg. 1 for 10 iterations provides a good trade-off between runtime and accuracy, and that more iterations lead to comparable bounds. Assuming a fixed amount of iterations for the eigendecomposition, finding $\Delta$ for DS* and $\lambda_{\min}^\star$ for DS++ have equal asymptotic complexities, and thus the scalability of both is comparable. Nevertheless, DS++ is generally faster (Table 2), whereas DS* achieves tighter bounds (Fig. 4).

In order to efficiently find a $\Delta$ that leads to a good convex relaxation, in Alg. 1 we fixed $d{=}\mathbf{0}$ and optimised over $d_1$ and $d_2$. While Alg. 1 can easily be extended to also find $d$, our preliminary experiments with such an approach led to slightly better lower bounds, but to considerably worse upper bounds. We leave an in-depth exploration of using full matrices $D_1$ and $D_2$ as well as $d{\neq}\mathbf{0}$ for future research.

## 10. Conclusion

We have presented a general convex relaxation framework for quadratic optimisation problems over permutations. In contrast to lifting-based convex relaxation methods that use variables of dimension $\mathcal{O}(n^4)$, our approach does not increase the number of variables to obtain a convex relaxation and thus works with variables of dimension $\mathcal{O}(n^2)$. Moreover, our approach is at least as tight as existing (lifting-free) relaxation methods as they are contained as special cases. To achieve our relaxation we have analysed a class of parametrised objective functions that are equal over permutation matrices, and we provided insights on how to obtain parametrisations that lead to tighter convex relaxations. In particular, we have introduced a proximal subgradient type method that is able to efficiently approximate such a parametrisation. Overall, we have presented a powerful framework that offers a great potential for future work on convex and non-convex methods for diverse matching problems, which is confirmed by our experimental results.